\documentstyle[12pt]{article}

\makeatletter
\def\eqalign#1{\null\vcenter{\def\\{\cr}\openup\jot\m@th
  \ialign{\strut$\displaystyle{##}$\hfil&$\displaystyle{{}##}$\hfil
      \crcr#1\crcr}}\,}
\makeatother

\newcommand{\be}{\begin{equation}} 
\newcommand{\ee}{\end{equation}}
\newcommand{\beq}{\begin{eqnarray}}
\newcommand{\eeq}{\end{eqnarray}}
\newcommand{\bt}{\begin{theorem}}
\newcommand{\et}{\end{theorem}}
\newcommand{\bl}{\begin{lemma}}
\newcommand{\el}{\end{lemma}}
\newcommand{\bc}{\begin{corollary}}
\newcommand{\ec}{\end{corollary}}
\newcommand{\bp}{\begin{prop}}
\newcommand{\ep}{\end{prop}}
\newcommand{\ba}{\begin{array}}
\newcommand{\ea}{\end{array}}

\newcommand{\la}{\label}
\newcommand{\ci}{\cite}

\newcommand \qed {\hskip 6pt\vrule height6pt width5pt depth1pt \bigskip}
\newtheorem{theorem}{THEOREM}
\newtheorem{lemma}[theorem]{LEMMA}
\newtheorem{corollary}[theorem]{COROLLARY}
\newtheorem{prop}[theorem]{PROPOSITION}

\newcommand{\lb}{\lambda}

\newcommand{\th}{\theta}

\newcommand{\bi}{\bibitem}

\newfont{\msbm}{msbm10 scaled\magstep1}
\newfont{\msbms}{msbm7 scaled\magstep1} 

\newcommand{\bbr}{\mbox{$\mbox{\msbm R}$}}

\newcommand{\bbz}{\mbox{$\mbox{\msbm Z}$}}

\begin{document}
\today
\bigskip\bigskip\bigskip
\begin{center}
{\Large\bf Continuity of the measure of the spectrum for discrete 
quasiperiodic operators}\\
\bigskip\bigskip\bigskip\bigskip
{\large S. Ya. Jitomirskaya$^1$, I. V. Krasovsky$^2$}\\
\bigskip 
$^1$ Department of Mathematics, University of California\\
Irvine, CA 92697, USA\\
E-mail: szhitomi@math.uci.edu\\

$^2$ Technische Universit\"at Berlin,   
Institut f\"ur Mathematik MA 7-2\\
Strasse des 17. Juni 136, D-10623, Berlin, Germany\\
E-mail: ivk@math.tu-berlin.de\\
\bigskip\bigskip
 
\end{center}
\bigskip\bigskip\bigskip

\noindent
{\bf Abstract.}
We study discrete Schr\"odinger operators 
$(H_{\alpha,\theta}\psi)(n)=
\psi(n-1)+\psi(n+1)+f(\alpha n+\theta)\psi(n)$ on $l^2(Z)$, 
where $f(x)$ is a real analytic  periodic function of period 1.
We prove a general theorem relating the measure of the spectrum of 
$H_{\alpha,\theta}$
to the measures of the spectra of its canonical rational approximants 
under the condition that the Lyapunov exponents of $H_{\alpha,\theta}$
are positive. For the almost Mathieu operator ($f(x)=2\lambda\cos 2\pi x$)
it follows that the measure of the spectrum is equal to $4|1-|\lambda||$ for
all real $\theta$, $\lambda\ne\pm 1$, and all irrational $\alpha$.

\newpage
\section{Introduction}
Consider quasiperiodic  operators acting on $l^2(\bbz)$ and given by:
\begin{equation}
(H_{\alpha,\theta}\psi)(n)=
\psi(n-1)+\psi(n+1)+f(\alpha n+\theta)\psi(n),
\qquad n=\dots,-1,0,1,\dots,\label{1}
\end{equation}
where $f(x)$ is a real  analytic periodic function of period 1. 
Denote by $S(\alpha,\theta)$ the spectrum of $H_{\alpha,\theta}.$ 
For rational $\alpha=p/q$ the spectrum 
consists of at most $q$ intervals. Let $S(\alpha)= \bigcup_{\theta\in\bbr}
S\left(\alpha,\theta\right).$  Note that for irrational $\alpha$
the spectrum of $H$ (as a set) is independent of $\theta$ (see, e.g., 
\ci{Cycon}), and therefore $S(\alpha,\theta)=S(\alpha).$ In this paper 
we study continuity in $\alpha$ of $S(\alpha)$ and its
measure. For sets, we will use $|\cdot |$ to denote the Lebesgue measure. 

The fact that various quantities could be easier to analyse and sometimes 
are even computable
for periodic operators, $H_{p/q,\theta}$,  makes results on continuity in 
$\alpha$ particularly important. 
For example, the Aubry-Andre conjecture on the measure of the 
spectrum \ci{aa} states that for the almost Mathieu operator given by
(\ref{1}) with $f(\th)=2\lb\cos2\pi\th$, for irrational $\alpha$ and 
all real $\lb,\th$ there is an equality
\be
\la{aa}
|S_\lambda(\alpha,\th)|=4|1-|\lb||.
\ee
Avron,van Mouche, Simon \ci{ams} proved that, for $|\lb| \not= 1,$ 
$|S_\lambda(p_n/q_n)|\to 4|1-|\lb||$ as $q_n\to \infty$, and
Last \ci{l2} established this fact for $|\lb| = 1$. Given these theorems, 
the proof of the Aubry-Andre conjecture was reduced to a
continuity result.
   
   The continuity in $\alpha$ of $S(\alpha)$ was proven in \ci{as,el}.
Continuity of the measure of the spectrum is  a more delicate issue, since, 
in particular, $|S(\alpha)|$ can be discontinuous at
rational $\alpha.$ Establishing continuity at irrational $\alpha$ requires 
quantitative estimates on the continuity of the
spectrum. The first such result, namely the H\"older-${1\over 3}$ 
continuity was proved in
\ci{ecy}. That argument was improved to the H\"older-1/2 continuity
(for arbitrary $f \in C^1$) in
\ci{ams} and subsequently used in \ci{l1,l2} to establish (\ref{aa})
 for the almost Mathieu operator 
for $\alpha$ with unbounded continuous fraction expansion, therefore proving 
the Aubry-Andre conjecture for a.e. (but not all)
$\alpha.$ It was essentially argued in
\ci{ams} that H\"older continuity of any order larger than 1/2 would imply 
the desired continuity property of the measure of the
spectrum for all $\alpha$. Such continuity (more precisely, almost 
Lipschitz continuity, as in Theorem 3) was proved in \ci{jl}
for the almost Mathieu operator with $|\lb|>14.5$ (or dual regime), using 
exponential localization for Diophantine frequencies and
delicate analysis of generalized eigenfunctions 
(positions of resonances and behavior between them). 

In the present paper we prove almost Lipschitz continuity of the spectrum 
for operators (\ref{1}) with arbitrary analytic $f$
under the assumption of positivity of the Lyapunov exponents. We would like 
to note that, unlike in \ci{jl}, we do not rely on
the theorem on exponential localization, and the present paper is essentially 
self-contained. In fact, the localization theorem for
general analytic potentials, proved in
\ci{bg}, establishes pure point spectrum for a.e. $\alpha$ without explicit 
Diophantine control, and therefore would not be
applicable for our purposes. On the other hand, we would like to mention 
that, as the proof below shows, continuity of the spectrum
requires simpler arguments than localization.

Let $M_n$ be the transfer-matrix of
$H\psi=E\psi$ :
\[
{\psi(n+1)\choose\psi(n)}=M_n(\theta,E){\psi(n)\choose\psi(n-1)}.
\]
Then 
\[
 M_n(\theta,E)=\pmatrix{E-f(\alpha n+\theta)& -1\cr 1& 0},\qquad
 n=\dots,-1,0,1,\dots. 
\]
Let $T_n(\theta,E)=M_{n-1}(\theta,E)\cdots
M_1(\theta,E)M_0(\theta,E)$ be the $n$-step transfer-matrix.
The Lyapunov exponent $\gamma(E,\alpha)$ 
of the family (\ref{1}) is defined as follows:
\begin{equation}
\gamma(E,\alpha)\equiv\lim_{n\to\infty}{1\over n}
\int_0^1 \ln||T_n(\theta,E)||d\theta=
\inf_n{1\over n}\int_0^1 \ln||T_n(\theta,E)||d\theta.\label{gamma}
\end{equation}
It is well defined by (\ref{gamma}) by Kingman's subadditive ergodic theorem 
(see, e.g., \ci{Cycon}). Note that ${\mathrm det}\,M_n(\theta,E)=1$ 
whence it follows that $\gamma(E,\alpha)\ge 0$. 
Our main result is

{\bf Theorem 1.} {\it Let $\gamma(E,\alpha)>0$ for Lebesgue a.e. $E$. 
Then for any real $\theta$,
\begin{equation}
|S(\alpha,\theta)|=\lim_{n\to\infty}\left|
S\left(\frac{p_n}{q_n}\right)\right|,
\end{equation}
where ${p_n/q_n}$ is the sequence of canonical rational approximants 
to $\alpha$.}\\

{\bf Remarks.}
\begin{enumerate}
\item
Note that a.e. positivity of $\gamma$ implies that $\alpha$ is irrational.
\item For $\alpha$ whose coefficients of
the continued fraction expansion  form an unbounded sequence, the condition 
$\gamma(E,\alpha)>0$ is redundant.
The result of Theorem 1 follows in this case (by a simple 
generalization of the argument in \ci{l1}) from the H\"older-$1/2$ 
continuity of the spectrum established in \ci{ams}. For $\alpha$'s whose 
continued fraction coefficients form a bounded sequence (henceforth, we 
denote the set of such $\alpha$ by $\Omega$) this continuity is not 
sufficient to imply Theorem 1. For such $\alpha$ we prove the theorem by 
first establishing the almost-Lipschitz continuity of the spectrum for which, 
in turn, we use positivity of the Lyapunov exponent.
\item Note that it is not clear apriori that the limit above exists.
\item While our proof can be easily adapted for arbitrary $\alpha$ we, 
for the reason above, concentrate only on the case $\alpha \in \Omega$, 
which slightly simplifies certain arguments.
\end{enumerate}

We would like to add that recently a number of remarkable properties of 
quasiperiodic operators with analytic potentials have been
established based on  the positivity of the  Lyapunov exponents 
\ci{j,bg,jl2,gs}. Theorem 1, therefore, adds to the collection of
general corollaries of positive Lyapunov exponents for analytic potentials. 


For the almost Mathieu operator Theorem 1 immediately implies 

{\bf Theorem 2.} {\it For every $\lambda$, $\theta\in\bbr$,
$|\lambda|\ne 1$, and every irrational $\alpha$, the measure of the
spectrum of the almost Mathieu operator $H_{\lambda,\alpha,\theta}$
is equal to $4|1-|\lambda||$.}

This improves some results of \ci{l1,jl} and establishes the Aubry-Andre 
conjecture on the measure of the spectrum for all
values of parameters in the noncritical case. For the critical case, $\lambda =2,$
zero measure of the spectrum is known for the full measure set of $\alpha$ \ci{hs,l2},
however extending it to all irrational $\alpha$ remains an open problem.


{\it Proof:}
It was shown in \ci{ams} that on a sequence of 
rational approximants $p_n/q_n$ to any real $\alpha$ we have that 
$|S_\lambda( p_n/q_n)|$ 
($|\lambda|\ne 1$) converges to $4|1-|\lambda||$. Assume $\lambda>0$
as the case $\lambda<0$ is easily obtained by the transform 
$\theta\to\theta+1/2$ and the case $\lambda=0$ is trivial.
We further assume $\lambda>1$ as the case $0<\lambda<1$
can than be analyzed using duality (see, e.g. \ci{Cycon}).
For $\lambda>1$, it is known (e.g., \ci{Cycon})
that $\gamma(E,\alpha)>0$ for all real $E$ and 
irrational $\alpha$, which  concludes the proof in view of Theorem 1.\qed

In Section 2 we establish almost Lipschitz continuity of the spectrum and 
in Section 3 we complete the proof of Theorem 1.

\section{Continuity of the spectrum.}
In what follows we shall often omit the indices
$\alpha,\theta$, or $E$ from the notation. 
Denote the basis in which $H$ has the tridiagonal structure (\ref{1})
by $\{e_i\}_{i=-\infty}^\infty$. Let 
$(E-H)_{[x_1,x_2]}$ be the operator $E-H$ restricted to the
subspace spanned by $\{e_i\}_{i=x_1}^{x_2}$. Set $k=x_2+1-x_1$. Consider 
the polynomial in $E$:
\[
\widetilde P_k(\alpha,\theta
,E)=\det((E-H)_{[0,k-1]}),\qquad \widetilde P_0=1.\label{P}
\]
Let $f_k$ be the Fourier  truncation of  $f(\theta):$ 
\begin{equation}
f_k(\theta)=\sum_{j=-k^2}^{k^2}
\hat f_j e^{2\pi ij\theta}.\label{Fourier}
\end{equation}
As $f(\theta)$ is an analytic function, its Fourier coefficients $\hat f_j$
are of order $e^{-cj}$. Hence
$|f(\theta)-f_k(\theta)|<\exp(-dk^2)$, where $d>0$ depends 
on $f(\theta)$, for $k$ sufficiently large.
Set $z=\exp(2\pi i\theta)$.
Let $H_k$ be the operator $H$ with $f(\theta)$  replaced by 
$f_k(\theta)$ and consider
\begin{equation}
P_k(\alpha,z,E)=
z^{k^3}\det((E-H_k)_{[0,k-1]}).\label{P2}
\end{equation}
Obviously, $P_k(z)$ is a polynomial of degree $2k^3$
in $z$. Using induction for the three-term recurrence relation (which connects
$\widetilde P_n$, $\widetilde P_{n-1}$, and $\widetilde P_{n-2}$) 
and boundedness of $f(\theta)$,
we obtain for all $z$, $|z|=1$, and  $E$ in a bounded set:
\begin{equation}
|P_k|-e^{-d k^2}<|\widetilde P_k|<
|P_k|+e^{-d k^2},\qquad |P_k|<e^{Dk},\label{ineq}
\end{equation}
where $D<\infty$, $d$ is somewhat decreased from the prior value
but still positive, and $k$ is larger than some $k_1$. 

We shall now obtain the upper and lower bounds on $|P_k(z,E)|$ which we
use later to establish the continuity properties of the spectrum.

{\bf Lemma 1.} (Estimate of $|P_k(z,E)|$ from below) 
{\it For every $\alpha\in\Omega$, $E\in S,$
$\varepsilon'>0$ and at least one integer $k$ out of each three consecutive
integers $s$, $s+1$, $s+2$, $s>k_1,$ there exists 
$C(\varepsilon',\alpha)$ such that in any interval of length $Ck^3$ 
there is an integer $m$ such that
$|P_k(e^{2\pi i (\th + m\alpha)},E)|\ge e^{(\gamma(E)-\varepsilon')k}$. 
The constant $C$ is independent of $E$.}

{\bf Remark.} This statement, appropriately adjusted, holds under any 
Diophantine condition on $\alpha$.

{\it Proof:}
  Since, as is easy to verify, 
\[
T_k(\theta,E)= \pmatrix{
\widetilde P_k(\theta)& -\widetilde P_{k-1}(\theta+\alpha)\cr
\widetilde P_{k-1}(\theta)& -\widetilde P_{k-2}(\theta+\alpha)},\qquad k\ge 2,
\]
we can rewrite (\ref{gamma}) in the form 
\begin{equation}
\gamma(E)=
\inf_k{1\over k}\int_0^1 \ln\mu_k(\theta)d\theta,\label{gamma2}
\end{equation}
where $\mu_k(\theta)=\max\{|P_k(\theta)|,|P_{k-1}(\theta)|,
|P_{k-1}(\theta+\alpha)|,
|P_{k-2}(\theta+\alpha)|\}.$ 
Let $B_k=\{\theta\in[0,1):
 \mu_k(\theta)\ge e^{(\gamma(E)-\varepsilon')k}\}$. Then, using (\ref{ineq})
and (\ref{gamma2}), we obtain
\[
k\gamma(E)\le\int_0^1 \ln\mu_k(\theta)d\theta=\int_{B_k}+
\int_{[0,1)\setminus B_k}\le |B_k|kD+(1-|B_k|)(\gamma(E)-\varepsilon')k.
\]
Therefore, 
$|B_k|\ge\varepsilon'/(D-\gamma(E)+\varepsilon')\ge 3c_1(\varepsilon')$ 
for all E in any bounded set (in particular, the
spectrum of $H.$)  Hence among each three consecutive indices 
$s,s+1,s+2>k_1$, there is at least one, denote it $k$, such that 
\[
|N_k:=\{\theta\in[0,1):
 |P_k(e^{2\pi i\theta})|\ge e^{(\gamma(E)-\varepsilon')k}\}|\ge 
c_1(\varepsilon')>0
\]
Denote by $M$ the set of all such $k$.

Since $P_k$ is a polynomial of degree $2k^3$ in $e^{2\pi i\theta}$, the set
$N_k$ consists of less than $9k^3$ intervals. 
Let $J$ be the
interval in $N_k$ with the maximum length.
We have
\begin{equation}
|J|\ge\frac{c_1}{9k^3}.
\end{equation}
It follows from the Diophantine properties of $\Omega$ that  
$\theta+m\alpha({\mathrm mod}\,1)\in J$ for at least one $m\in I$,
where $I$ is any interval of length $Ck^3.$ For reader's convenience 
we give a simple argument for the above.

Let $\alpha\in[0,1)$ be an irrational and consider the
continued fraction expansion $\alpha=[0;a_1,a_2,\dots]$, where $a_i$
are positive integers. Recall (e.g.,\ci{Khintchine}) that the 
convergents $p_n/q_n=[0;a_1,a_2,\dots,a_n]$ have the properties
\begin{equation}
\alpha=\frac{p_n}{q_n}+\frac{\widetilde\delta}{q_n q_{n+1}},
\qquad |\widetilde\delta|<1,\qquad
q_n=a_nq_{n-1}+q_{n-2},\qquad n=3,4,\dots.\label{cf}
\end{equation}
and the greatest common divisor $gcd(p_n,q_n)=1$.

For a fixed $n$, consider the set $\th_j= \th+\alpha j({\mathrm mod}\,1)$, 
$j=0,1,\dots, q_n-1$. Let $\th_{j}$ and $\th_{j'}$ be nearest neighbors. 
We have, using (\ref{cf}),
\begin{equation}\la{md}
\left|\alpha j-\frac{jp_n}{q_n}\right|=
\left|\frac{j\widetilde\delta}{q_nq_{n+1}}\right|<
\frac{j}{q_nq_{n+1}}<\frac{1}{q_n}.
\end{equation}
Since the points
${jp_n\over q_n}({\mathrm mod}\,1)$, $j=0,1,\dots, q_n-1$ are the same as
${j\over q_n}({\mathrm mod}\,1)$, $j=0,1,\dots, q_n-1,$ up to the ordering,
 we have, using (\ref{md}):
\begin{equation}
\max|\theta_j-\theta_{j'}|< {3\over{q_n}},
\qquad j=1,2,\dots,q_n.
\end{equation}
Therefore, if $c_1/9k^3\ge 3/q_n,$ there exists $0\le m\le q_n-1$ such that
$\theta_m\in J$. Accordingly, choose $n$ so that
\begin{equation}
q_{n-1}<\frac{27 k^3}{c_1}\le q_n.\label{defqn}
\end{equation}

For $\alpha \in \Omega$ we have that $a_n\le B(\alpha)$ for all $n$. Then
\begin{equation}
q_n=a_nq_{n-1}+q_{n-2}<2 B(\alpha)q_{n-1}<C k^3.\label{ck3}
\end{equation}
where $C=60B(\alpha)/c_1$.
Therefore, at one of any $q_n$ (as defined by (\ref{defqn})) of
the points $\th_m$ with consecutive indices  we have the
needed estimate for $P_k(e^{2\pi i \th})$. Inequality (\ref{ck3}) implies 
the statement of the Lemma. \qed

{\bf Lemma 2.} (Estimate of $|P_k(z,E)|$ from above) 
{\it For any  $\delta,\varepsilon >0$,
there exists a set $F(\delta,\varepsilon)\subset S$ such 
that $|F(\delta,\varepsilon)|<\delta$ and for any $\alpha\in\Omega$,
$E\in S\setminus F(\delta,\varepsilon)$, and sufficiently 
large $k$ ($k>k_2(\delta,\alpha,\varepsilon)$) we have
$|P_k(z,E)|\le e^{(\gamma(E)+\varepsilon)k}$. The constant $k_2$ is
independent of $z$ and $E$.}

{\bf Remark.} Lemma 2 actually holds uniformly, for all $E$ in a bounded set, \ci{bj}.
However, this result requires a relatively 
complicated argument, while a nonuniform
statement above follows immediately from a general theorem on uniquely ergodic dynamical systems and is sufficient for our 
purposes. 
 

{\it Proof:}  
It is proved in
\ci{f} that for any continuous subadditive cocycle $f_n$ on a uniquely ergodic system, $\limsup \frac1n f_n (x) \le \lim 
\frac 1n \int
f_n d\mu(x)$  uniformly in $x.$ Therefore, for large $k$
we have $|P_k(z)|<e^{(\gamma_k(E,\alpha)+\varepsilon)k}$ for all $|z|=1.$ 

Let $D_n=\{E\in S: \forall k>n\ \forall |z|=1,\ |P_k(z,E)|\le 
e^{(\gamma(E)+\varepsilon)k}\}$. Then
we have $|S\setminus D_n|\to 0$ as
$n\to\infty$, which completes the proof of Lemma 2.\qed

We are now ready to formulate the result about continuity of the
spectrum.

{\bf Theorem 3.}
{\it Suppose $\alpha\in\Omega$ and $\gamma(E,\alpha)>0$ for Lebesgue a.e. $E$.
Then for any  $\delta>0$,
there exists a set $A(\delta)\subset S(\alpha)$ such that 
$|A(\delta)|<\delta$ and for any 
$E\in S(\alpha)\setminus A(\delta)$ and sufficiently small 
difference $|\alpha-\alpha'|<h(\alpha,\delta)$ we can find
$E'\in S(\alpha')$
 such that
\begin{equation}
|E'-E|< c(\alpha,\delta)|(\alpha-\alpha')\ln^3|\alpha-\alpha'||.\label{cont}
\end{equation}
The constant $c$ is independent of $E$, $E'$, $\alpha'$.}\\

{\bf Remarks.}\begin{enumerate}\item The exponent $3$ in (\ref{cont}) 
can be replaced by $2+\epsilon$ with
any $\epsilon>0$, as the reader can easily verify starting with equation 
(\ref{Fourier}). For the almost Mathieu case, $3$ can be replaced by $1$.
\item A similar statement (with $3$ replaced by a higher power) holds under 
a standard Diophantine condition on $\alpha.$
\item This statement can be made uniform in $E\in S$ if we require 
$\gamma(E)$ to be bounded away from 0 (as in the case of the
almost Mathieu operator) and use the uniform version of Lemma 2.
\end{enumerate}

 {\it Proof:} As in \ci{ecy,ams,l1,jl}, for a given 
$E \in S(\alpha)$ we construct an approximate eigenvector for
$H_{\alpha}.$  In order to obtain almost-Lipschitz continuity we need the 
error of approximation to be exponentially small, as in \ci{jl}.

Fix $\alpha\in\Omega$. Let $G=\{E\in S:\gamma(E)<g\}$
for some $g>0$. Then $|G|\to 0$ as $g\to 0$.
Let $g$ be such that $|G|<\delta/2$. Choose positive $\varepsilon'$, 
$\varepsilon''$ so that $\varepsilon'+\varepsilon''\le g/16$.
Set $A(\delta)=F({\delta\over 2},\varepsilon'')\cup G$.
Fix $E\in S\setminus A(\delta)$.
Set $K=\{k\in M, k>2k_2(\delta,\alpha,\varepsilon'')+3\}$
with $k_2$ from Lemma 2, and $M$ defined in the proof of Lemma 1.

For any $E_0\in S$ 
let $G(x,y)$ be the matrix elements of $G=((H-E_0)_{[x_1,x_2]})^{-1}.$ 
Using Cramer's rule and (\ref{ineq}) at the last step,
we obtain:
\begin{equation}
|G(x,x_1)|<\frac{|P_{x_2-x}(e^{2\pi i\th_{x+1}},E_0)|+e^{-dk^2}}
{|P_k(e^{2\pi i\th_{x_1}},E_0)|-e^{-dk^2}},\qquad
|G(x,x_2)|<
\frac{|P_{x-x_1}(e^{2\pi i\th_{x_1}},E_0)|+e^{-dk^2}}
{|P_k(e^{2\pi i\th_{x_1}},E_0)|-e^{-dk^2}},\label{G}
\end{equation}
where $k=x_2-x_1+1>k_1$.

We shall now fix $x_1$, $k$, and $E_0$.
First take a $k\in K$. By Lemma 2, 
$|P_{[(k+j)/2]}(z,E)|<\exp{(\gamma(E)+\varepsilon'')(k/2+1)}$ for
$j=-3,-2,\dots,2,$ all $z$. We can find a neighborhood of $E,$ $r_k(E)$ of 
diameter smaller than 
$\exp(-kg/4),$ so that for any $E''\in r_k(E)$ we have for all
$z$ ($|z|$=1):
$|P_{[(k+j)/2]}(E'')|<|P_{[(k+j)/2]}(E)|\exp{(\varepsilon''(k/2+1))}$ and
$|P_k(E'')|>|P_k(E)|\exp(-\varepsilon'k)$.
Now take a generalized eigenvalue $E_0\in r_k(E)$. Let $\psi$ be the
corresponding generalized eigenfunction, that is a formal
solution to the equation $H\psi=E_0\psi$ with $|\psi(x)|<\widetilde C (|x|+1)$ 
for all $x.$  By \ci{Shnol,Cycon} the set of $E_0$ which admit such
solutions 
is dense in the spectrum. 
There holds:
\begin{equation}
\max_x\frac{|\psi(x)|}{|x|+1}=\frac{|\psi(x_{\max})|}{|x_{\max}|+1}=R<\infty
\end{equation}
for some point $x_{\max}$. We normalize $\psi$ so that $R=1$.

As is easily seen by considering
$(H-E_0)_{[x_1,x_2]}\psi$, we have for $x\in[x_1,x_2]$:
\begin{equation}
-\psi(x)=G(x,x_1)\psi(x_1-1)+G(x,x_2)\psi(x_2+1).\label{psi}
\end{equation}
Define the interval 
$I=[x_{\max}-k-Ck^3,x_{\max}-k]$, where $C=C(\varepsilon', \alpha)$ 
is the constant from Lemma 1. Then $|I|=Ck^3$, and in view of 
Lemma 1, we can find $x_1=m$ in $I$ (which fixes
position of the interval $[x_1,x_2]$) such that 
$|P_k(e^{2\pi i\th_m},E)|\ge\exp{(\gamma(E)-\varepsilon')k}$ 
and hence
$|P_k(e^{2\pi i\th_m},E_0)|>\exp{(\gamma(E)-2\varepsilon')k}$. 
We now evaluate $\psi(x)$ at the midpoint of the interval
$[x_1,x_2]$: $x_0=x_1+[(k-1)/2]$ and at its nearest neighbors.
Using Lemma 1 to evaluate the denominator 
in (\ref{G}), we obtain
$|G(x_0,x_1)|<2\exp{((-\gamma(E)+\varepsilon)k/2+\gamma(E)+\varepsilon)}$, 
$\varepsilon=4(\varepsilon'+\varepsilon'')$
and the same estimate for $G(x_0,x_2)$ for sufficiently large 
$k \in K$ (depending on $d$).
 Applying (\ref{psi}), we get
\begin{equation}
|\psi(x_0-1)|,
|\psi(x_0)|< 4e^{-(g-\varepsilon)k/2+D+\varepsilon}(|x_{\max}|+(C+2)k^3)
\label{est}
\end{equation}
for sufficiently large $k\in K$.

Now let $I'=[x_{\max}+k,x_{\max}+k+Ck^3].$ Similarly, we
choose an interval $[x^{'}_1,x^{'}_2]$ so that 
at $x^{'}_1\in I'$ we have a lower bound of Lemma 1. We apply again Lemmas 1 and 2 to get the same estimate
(\ref{est}) for $\psi(x^{'}_0)$ and $\psi(x^{'}_0+1)$ 
at the midpoint $x^{'}_0$ of $[x^{'}_1,x^{'}_2]$ and at $x^{'}_0+1$.
By construction, $L=|x_0-x^{'}_0|$ satisfies $k<L<2(C+1)k^3$.
Set
\[
\psi_L(x)=\cases{\psi(x), & $x\in[x_0,x^{'}_0]$\cr 0, & otherwise}.
\]
and $\phi_L(x)=\psi_L(x)/||\psi_L||$.
Since $x_{\max}\in [x_0,x^{'}_0]$, we have $||\psi_L||>|x_{\max}|+1$.
Hence
\begin{equation}
\eqalign{
\frac{|\psi(x_0-1)|}{||\psi_L||},
|\phi_L(x_0)|,
\frac{|\psi(x^{'}_0+1)|}{||\psi_L||}, |\phi_L(x^{'}_0)|
< 4e^{-(g-\varepsilon)k/2+D+\varepsilon}
\frac{|x_{\max}|+(C+2)k^3}{|x_{\max}|+1}<\\
(C+2)e^{-(g-2\varepsilon)k/2}}
\label{est2}
\end{equation}
for $k\in K$, $k>K_0$, where $K_0(\varepsilon,D,d)$ is 
sufficiently large.

By the variational principle, there exists a point $E'$ in the
spectrum of $H_{\alpha',\theta'}$ (here $\theta'=
(\alpha-\alpha') x_{\max}+\theta$)  such that
\begin{equation}
\eqalign{
|E'-E|\le||(H_{\alpha',\theta'}-E)\phi_L||\le\\
||(H_{\alpha',\theta'}-H_{\alpha,\theta})\phi_L||+
||(H_{\alpha,\theta}-E_0)\phi_L||+|E-E_0|<\\
C'|\alpha-\alpha'|L+4(C+2)e^{-(g-2\varepsilon)k/2}
+e^{-kg/4}<\\
C' 2(C+1)|\alpha-\alpha'|k^3+4(C+3)e^{-gk/4},}\label{E}
\end{equation}
where we applied the estimate $|f(\alpha n+\theta)-f(\alpha' n+\theta')|<
C'|\alpha-\alpha'|L$ with some $C'>0$ and
used that, by our choice of parameters, $2\varepsilon\le g/2$. 
Let $\alpha'$ be sufficiently close to $\alpha,$ so that
$k=[4|\ln|\alpha-\alpha'||/g]$ (or $k\pm 1)$ is in $K$ and larger than $K_0$
and $|\ln|\alpha-\alpha'||$ is sufficiently large depending on values of 
$C$, $C'$, and $g$.
Then we obtain from (\ref{E}) the statement of Theorem 3 with 
$c(\alpha,\delta)=2^7C'(C+1)/g^3+1$.\qed

\section{Proof of Theorem 1}
As we noted in the introduction,
the theorem for $\alpha\notin\Omega$ is easy to prove following 
\ci{ams,l1}.
Take $\alpha\in\Omega$ and consider 
the sequence of its canonical rational approximants $p_n/q_n.$ 
Because of continuity in $\theta$,
the set $S(p_n/q_n)$ consists of at most 
$q_n$ disjoint intervals, say
$S(p_n/q_n)=
\cup_{i=1}^{{q'}_n}[a_i,b_i]$, ${q'}_n\le q_n$.
Given the continuity result of Theorem 3, the proof of Theorem 1 
differs from the corresponding analysis  
for the almost Mathieu operator in \ci{l1} only in one detail. 
By \ci{T83,l2} we have   
$|S(\alpha)|\ge\limsup_{n\to\infty}|
S(p_n/q_n)|$. We are going to show that $|S(\alpha)|\le
\liminf_{n\to\infty}|S(p_n/q_n)|$
when $\gamma(E)>0$ for a.e. $E$. 
For all $n$ sufficiently large (such that $|p_n/q_n-\alpha|<h(\alpha,\delta)$),
Theorem 3 says that, except perhaps for the points belonging to
$A(\delta)$, all the other points in the spectrum of $H_\alpha$ are
necessarily inside the set ($\alpha'=p_n/q_n$)
\[
\cup_{i=1}^{{q'}_n}[a_i-
c|(\alpha-\alpha')\ln^3|\alpha-\alpha'||,
b_i+c|(\alpha-\alpha')\ln^3|\alpha-\alpha'||].
\]
Therefore,
\begin{equation}
|S(\alpha)|<
|S(p_n/q_n)|+
2 c q_n|(\alpha-p_n/q_n)\ln^3|\alpha-p_n/q_n||+\delta
\end{equation}
As $n\to\infty$,
the second addend at the right hand side tends to 0 for any irrational 
$\alpha$ because of (\ref{cf}). The constant $\delta>0$ can be taken arbitrary 
small  by Theorem 3. Hence, 
$|S(\alpha)|=\lim_{n\to\infty}|S(p_n/q_n)|.$\qed

\section{Acknowledgements}
We would like to acknowledge the support of the  NSF grant DMS-9800860. 
S.J. was also supported under NSF grant DMS-0070755, and I.K. by 
Sonderforschungsbereich 288.
I.K. is grateful to S. Jitomirskaya and A. Klein for their hospitality in the
University of California, Irvine, where part of this work was written.

\end{document}